\documentclass[letterpaper,reqno]{amsart}

\usepackage{amsmath} 
\usepackage{amsfonts} 
\usepackage{amssymb} 
\usepackage{epsfig} 
\usepackage{graphicx} 
\usepackage{multirow} 
\usepackage{verbatim} 
\usepackage{rotating} 
\usepackage[all]{xy} 
\usepackage{appendix}
\usepackage[english]{babel}
\usepackage{epigraph}
\usepackage{mathtools}

\newtheorem{theorem}{Theorem}[section]
\newtheorem{lemma}[theorem]{Lemma}

\theoremstyle{definition}
\newtheorem*{theorem*}{Theorem}
\newtheorem{definition}[theorem]{Definition}
\newtheorem{proposition}[theorem]{Proposition}

\newtheorem{warning}[theorem]{Warning}

\theoremstyle{remark}
\newtheorem{remark}[theorem]{Remark}
\theoremstyle{notation}
\newtheorem{notation}[theorem]{Notation}
\numberwithin{equation}{section}
\theoremstyle{corollary}

\newtheorem{conjecture}[theorem]{Conjecture}
\usepackage[bookmarks=false]{hyperref}
\newcommand{\Map}{\mathrm{Map}}

\newcommand{\Top}{\mathsf{Top}}

\newcommand{\Gr}{\mathsf{TGr}}
\newcommand{\Homology}{\mathrm{H}}
\newcommand{\Ho}{\mathrm{Ho}}
\newcommand{\E}{\mathbf{E}}
\newcommand{\Lc}{\mathbf{L}}
\newcommand{\U}{\mathbf{U}}

\newcommand{\B}{\mathbf{B}}
\newcommand{\C}{\mathbf{C}}
\newcommand{\F}{\mathbf{F}}
\newcommand{\M}{\mathbf{M}}
\newcommand{\Path}{\mathbf{P}}

\newcommand{\Grp}{\mathsf{Gr}}

\newcommand{\hocolim}{\mathrm{hocolim}}

\newcommand{\colim}{\mathrm{colim}}

\newcommand{\Sing}{\mathbf{Sing}}

\begin{document}

\title[]{Friedlander-Milnor Conjectures}

\author[]{Ilias Amrani}

\address{Saint Petersburg State University, Russian Federation.}
\curraddr{Faculty of Mathematics and Mechanics. }
\address{Academic University of Saint Petersburg, Russian Federation.}
\curraddr{Department of Mathematics and Information Technology.  }
\email{i.amrani@spbu.ru}
\email{ilias.amranifedotov@gmail.com}






\begin{abstract}
We prove the Milnor conjecture for Lie groups and the Friedlander conjecture for complex algebraic Lie groups.   
 \end{abstract}
 
\maketitle
\tableofcontents
\section{Introduction}
\begin{conjecture}[Milnor \cite{milnor1983homology}]\label{conj1}
Let $G$ be a Lie group, and let $G^{\delta}$ be the underlying group (with discrete topology). The natural homorphism of topological groups $G^{\delta}\rightarrow G$ induces an isomorphism from the singular cohomology of the classifying space of $G$ to the singular cohomology of the classifying space of $G^{\delta}$ 
$$ \Homology^{\ast}(\B G;\mathbf{F}_{p})\rightarrow \Homology^{\ast}(\B G^{\delta};\mathbf{F}_{p}) $$ for any prime $p$. 
\end{conjecture}
There is also an other variation of the formulation of the previous conjecture in the algebraic setting, called the Friedlander conjecture.
\begin{conjecture}[Friedlander \cite{morel2010friedlander}]\label{conj2}
For any algebraically closed field $\mathbf{k}$, any reductive algebraic $\mathbf{k}$-group $G$ and any prime $p$ different from $char(\mathbf{k})$, the natural ring homomorphism
$$\Homology_{\textit{\'et}}^{\ast}(\B G;\mathbf{Z}/p\mathbf{Z})\rightarrow \Homology^{\ast}(\B G(\mathbf{k});\mathbf{Z}/p\mathbf{Z}) $$
from the $\mathbf{Z}/p\mathbf{Z}$-\'etale cohomology of the simplicial classifying space of G to the group cohomology with $\mathbf{F}_{p}=\mathbf{Z}/p\mathbf{Z}$ coefficients of the discrete group $G(\mathbf{k})$ is an isomorphism.
\end{conjecture}
\begin{remark}\label{remark0}
Using Grothendieck's comparison theorems between \'etale and singular cohomology \cite{grothendieck1972topos}, these two conjectures overlap when the group $G$ is a complex algebraic Lie group. For more details on the status of the conjectures we refer to \cite[Chapter V]{knudson2001homology} and \cite{morel2010friedlander}.
\end{remark}
\begin{remark}
The Friedlander conjecture was proven in \cite{friedlander1984cohomology} when $\mathbf{k}=\overline{\mathbf{F}}_{p}$ the algebraic closure of the finite field $\mathbf{F}_{p}$ of characteristic $p$.
\end{remark}

\subsection{Organization of the article and results}
In the section \ref{section1} and \ref{section3} we prove that the model category of topological groups is left proper. In section \ref{section2}, we introduce two functors $\Lc$ and $\E$. The hart of the article consists in two sections \ref{section5} and \ref{section6}, where we prove the following theorems:

\begin{theorem}\ref{crucial+}
Let $G$ be an object in $R^{\delta}/\Gr$ such that $\Lc(G)$ is a retract of a CW-complex as a space, then the natural map $$\mathrm{H}_{\ast}(\B G;k)\rightarrow \mathrm{H}_{\ast}(\B \Lc G;k)$$ is an isomorphism for any finite field $k$. 
\end{theorem}
\begin{theorem}\ref{pushout} 
Suppose that $G$ is a simple Lie group and $G^{\delta}$ is the same group with discrete topology then $G$ is isomorphic to $\Lc(G^{\delta})$ as topological group.
\end{theorem}
 
In the section \ref{section7}, we finally prove the (stable and unstable) Friedlander-Milnor conjectures:
\begin{theorem}\ref{main}
The Milnor conjecture is true, consequently the Friedlander conjecture is true for complex algebraic Lie groups.
\end{theorem}

\begin{theorem}(Stable Milnor conjecture \ref{stable})
 Let $p$ be any prime number and $\mathbf{k}$ the filed of real numbers $\mathbf{R}$ or the field of complex numbers $\mathbf{C}$. We define $\B GL_{\infty}(\mathbf{k})$ as $ colim_{n}\B GL_{n}(\mathbf{k})$, then the evident map $$\Homology_{\ast}(\B GL_{\infty}(\mathbf{k})^{\delta};\mathbf{F}_{p})\rightarrow \Homology_{\ast}(\B GL_{\infty}(\mathbf{k});\mathbf{F}_{p})$$ is an isomorphism. The same result holds if we replace $GL_{\infty}$ by the special linear group $SL_{\infty}$, by the orthogonal group $O_{\infty}$, by the special orthogonal group $SO_{\infty}$, by the unitary group $U_{\infty}$, by the special unitary group $SU_{\infty}$ or by the symplectic group $Sp_{\infty}$.
\end{theorem}

\section{The homotopy category of topological groups and monoids }\label{section1}
\def\cocartesien{%
    \ar@{-}[]+L+<-6pt,+1pt>;[]+LU+<-6pt,+6pt>%
    \ar@{-}[]+U+<-1pt,+6pt>;[]+LU+<-6pt,+6pt>%
  }

\def\cartesien{%
    \ar@{-}[]+R+<-1pt,+1pt>;[]+RD+<-1pt,+1pt>%
    \ar@{-}[]+D+<-1pt,+1pt>;[]+RD+<-1pt,+1pt>%
  }
\def\hcocartesien{%
    \ar@{-}^{h}[]+L+<-6pt,+1pt>;[]+LU+<-6pt,+6pt>%
    \ar@{-}[]+U+<-1pt,+6pt>;[]+LU+<-6pt,+6pt>%
  }
 \def\hcartesien{%
    \ar@{-}^{h}[]+R+<-1pt,+1pt>;[]+RD+<-1pt,+1pt>%
    \ar@{-}[]+D+<-1pt,+1pt>;[]+RD+<-1pt,+1pt>%
  }
\begin{warning}
A compactly generated space in our sense is the same as the notion of compactly generated space in the sense of \cite{strickland2009category}. While a compactly generated space in our sense is called $k$-space in the sense of Hovey \cite{hovey1999model}.
\end{warning}
The category of compactly generated (\textbf{CG}) topological spaces $\Top$ is a symmetric monoidal cofibrantly generated model category \cite{hovey1999model}. The homotopy category of topological spaces is denoted by $\Ho\Top$. The category of topological groups $\Gr$ is defined as the category of group objects in $\Top$. 
As for the category of sets and the category of groups we have an adjunction in the topological setting 
\begin{equation}
\xymatrix{\Top\ar@<2pt>[r]^-{\F} & \Gr \ar@<2pt>[l]^-{\U} }
\label{adj1} 
\end{equation}
We have chosen working with the category of compactly generated space (\textbf{CG}) for the following reason: First it guarantees that the colimits in the category of topological groups are the same as the colimits of underlying groups, once we apply the forgetful functor from the category of topological groups to the category of groups. Second the category of compactly generated spaces is monoidal closed model category. For more details about compactly generated spaces we refer to \cite{strickland2009category}. There is a small disadvantage working with compactly generated spaces which is the fact the cartesian product is a little bit modified (from topological viewpoint), hence our notion of topological groups differs slightly form the classical definition of topological groups. 
\begin{lemma}\label{model}
The category $\Gr$ is a cofibrantly generated model category where 
$G\rightarrow H$ is an weak equivalence (respectively a fibration) if and only if $\U G\rightarrow \U H$ is a weak equivalence (respectively a fibration) in the model category $\Top$. 
\end{lemma}
\begin{proof}
We apply the transfer principle developed in \cite{berger2003axiomatic} to the adjunction \ref{adj1}. The forgetful functor $\U$ preserves filtered colimits and since any object in $\Gr$ is fibrant it is enough to find a functorial path object in the category $\Gr$. 
For example, we take the functorial path object $$\Path (G)=\Map_{\Top}([0,1],G)$$  
equipped with compact-open topology which has a natural structure of topological group satisfying all the properties of a functorial path object.   
\end{proof}
\begin{remark}
The big advantage from considering the category of group objects in the category of compactly generated topological spaces is that the functor forgetting the topology $\mathbf{Forget}:\Gr\rightarrow \mathsf{Gr}$ from the category of topological groups to the category of groups has a right adjoint given by $\mathbf{Trivial}:\mathsf{Gr}\rightarrow \Gr$ sending a group to a topological group with trivial topology (the undiscrete topology). Therefore the functor $\mathbf{Forget}$ commutes with colimits. 
\end{remark}

\begin{remark}\label{grmon}
The model category $\Gr$ is a \textit{topological model category} in the sense of \cite[Definition 4.2, VII ]{EKMM}. It means that the model category $\Gr$ is enriched over the category of topological spaces $\Top$ and this structure is compatible with the model structures on $\Top$ and $\Gr$. For more details we refer to \cite[Theorem 4.4, VII ]{EKMM}. The model category of topological monoids $\mathsf{TMon}$ is also a a topological model category. Weak equivalences and fibrations in $\mathsf{TMon}$ are weak equivalences and fibrations in the underlying model category $\Top$. 
\end{remark}
For any topological groups $G$ and $H$, as a set, the mapping space $\Map_{\Gr}(G,H)$ is just the set of homomorphism in the category of topological groups but as a space it is the subspace of the topological space $\Map_{\Top}(\U G,\U H)$. This mapping space has a natural homotopical interpretation is the following way $\pi_{0}\Map_{\Gr}(G,H)\cong \Ho\Gr(G,H)$ if $G$ is a cofibrant topological group. The same remark holds for the topological model category $\mathsf{TMon}$ of topological monoids i.e. $\pi_{0}\Map_{\mathsf{TMon}}(M,N)\cong \Ho\mathsf{TMon}(M,N)$ if $M$ is a cofibrant topological monoid.\\
The category $\mathsf{sGr}$ of simplicial groups is a cofibrantly generated model category where the weak equivalences and fibration are weak equivalences and fibration in the underlying model category of simplicial sets (with the standard model structure).
\begin{lemma}\label{simplicial groups}
The following adjunction (geometric realization-singular functor)
\begin{equation}
\xymatrix{\mathsf{sGr}\ar@<2pt>[r]^-{|-|} & \Gr \ar@<2pt>[l]^-{\Sing} }
\label{adj11} 
\end{equation}
is a Quillen adjunction and a Quillen equivalence. 
\end{lemma}
\begin{proof}
It is a consequence of the fact that the category of simplicial sets and the category of topological spaces are Quillen equivalent
\cite[theorem 3.6.7]{hovey1999model}, and the fact that the functors $|-|$ and $\Sing$ are monoidal functors.
\end{proof}
\begin{remark}\label{remarkclass}
For any simplicial group $G_{\bullet}$, we define the classifying simplicial set by $\B_{\bullet}G_{\bullet}$, in \cite[Chapter V, section 4]{goerss2009simplicial}, the authors denote the functor $\B_{\bullet}$ by $\overline{\mathbf{W}}$. 
For any topological group $G$, we define the classifying space $\B G$ as a $|\B_{\bullet}\Sing G|$ where $\B_{\bullet}\Sing G$ is the classifying simplicial set for the simplicial group $\Sing G$. We can also define the classifying space of a simplicial group $G_{\bullet}$ as the diagonal of the bisimplicial set given by $[n][m]\mapsto \mathbf{N}_{n}G_{m}$ where $\mathbf{N}_{\bullet}$ is the simplicial nerve of a category and $G_{m}$ is seen as a category with one object $o$ and $Hom(o,o)=G_{m}$. In \cite{seb2008}, it is proven the two definitions of simplicial classifying spaces for simplicial groups are homotopy equivalent.    
\end{remark}
\begin{theorem}\label{classifying}
Let $\Top_{\ast}$ be the model category of topological pointed spaces (where weak equivalences and fibrations are those of the underlying model category $\Top$) and let $\Top^{0}_{\ast}$ be the category of pointed connected topological spaces,
the classifying functor $\B: \Gr\rightarrow \Top$ induces an equivalence of homotopy categories i.e., 
   $$\mathbf{B}: \Ho\Gr\rightarrow \Ho\Top^{0}_{\ast}$$ is an equivalence of categories. Moreover, for any (cofibrant) topological monoid $M$ and any topological group 
$G$ the natural map $$\Map_{\mathsf{TMon}}(M,G)\simeq \Map_{\Top_{\ast}}(\B M,\B G)$$
is a weak equivalence of topological spaces. Similarly for any (cofibrant) topological group $H$ 
$$\Map_{\mathsf{TGr}}(H,G)\simeq \Map_{\Top_{\ast}}(\B H,\B G)$$
is a weak equivalence of topological spaces.
\end{theorem}
\begin{proof}
It is a consequence of \cite[Proposition 1.11, proposition 1.12]{vogt2012homotopy} and \ref{simplicial groups}.
\end{proof}

\subsection{Cofibrations in the category $\Gr$}
The class of cofibrations has a particular issue in the model category $\Gr$ described in \ref{model}. In this section we will describe some of them. We recall that the generating set of cofibrations in the model category $\Top$ is given by $\{\partial \Delta^{n}\rightarrow \Delta^{n}, n\in \mathbf{N} \}$, in particular the set of generating cofibrations in $\Gr$ is given by  $$\{\F(\partial \Delta^{n})\rightarrow \F(\Delta^{n}), n\in \mathbf{N} \},$$ where $\Delta^{n}$ and $\partial \Delta^{n}$ are respectively the standard $n$-dimensional topological simplex and the boundary of the $n$-dimensional topological simplex. 
\begin{definition}\label{cellcof}
A \textbf{Cell-cofibration} $i:A\rightarrow B$ in the category $\Gr$ is a map form $A$ to a directed colimit indexed by a set $I$, more precisely 
 $A=A_{0}\rightarrow A_{1}\rightarrow \dots A_{s}\rightarrow A_{s+1} \rightarrow \dots$ where 
\begin{itemize}
\item $\colim_{I} (A_{0}\rightarrow A_{1}\rightarrow \dots A_{s}\rightarrow A_{s+1} \rightarrow \dots) = B$
\item the map $i$ is equal to the transfinite composition of the diagram bellow.
\item for any $s\in I$ the map  $A_{s}\rightarrow A_{s+1}$ is given by the pushout diagram of the form:
$$\xymatrix{\F(\partial\Delta^{n})\ar[r]\ar[d] & A_{s}\ar[d]\\
\F(\Delta^{n})\ar[r] & A_{s+1} \cocartesien
  } $$ 
for some $n\in \mathbf{N}$. 

\end{itemize}

\end{definition}
Any cofibration in the category $\Gr$ is a retract of a cell cofibration. 
\begin{remark}\label{laremarque}
Since colimits in the category $\Gr$ are colimits in the category of groups $\Grp$ if we forget about the topology, thus as a group, $A_{s+1}$ is isomorphic to the free product $A_{s}\ast \F(\Delta^{n}-\partial \Delta^{n})$. 
\end{remark}

\begin{lemma}\label{cofmon}
If $M\rightarrow N$ is a cofibration of topological monoids and as a space $M$ is a retract of a a $CW$-complex, then it is a cofibration of the underlying topological spaces.
\end{lemma}
\begin{proof}
First, recall that generating cofibrations in the category of topological monoids is given by the set $\{\M(\partial \Delta^{n})\rightarrow \M(\Delta^{n}): n\in \mathbf{N}\}$. Let $X\rightarrow Y$ be a cofibration of CW-complexs, then $\M(X)\rightarrow \M(Y)$ is a cofibration of topological monoids and cofibration of underlying spaces. Let $\M(X)\rightarrow M$ be any homomorphism of topological monoids, such that $M$ is a retract of a CW-complex. We consider the following pushout diagram in the category of topological monoids:

$$\xymatrix{\M(X)\ar[r]\ar[d] & M\ar[d] \\ 
\M(Y) \ar[r] & \M(Y)\ast_{\M(X)}M} $$  
in the proof of theorem 3.1 of \cite{batanin2017homotopy}, the authors prove that the map $s: M\rightarrow \M(Y)\ast_{\M(X)}M$ is a transfinite composition of cofibrations of spaces. For simplicity, we will call such map $s$ a primitive cell cofibration. Since any cell cofibration $M\rightarrow B$ in the category of topological monoids is a transfinite composition primitive cell cofibrations, it follows that $M\rightarrow B$ is a cofibration of underlying spaces. Finally, any cofibration $M\rightarrow A$ of topological monoids is a retract of a cell cofibration, therefore  $M\rightarrow A$ is a cofibration of underlying spaces. 
\end{proof}
\begin{proposition}\label{coftgroup}
If $G\rightarrow H$ is a cofibration of topological groups and as a space $G$ is a retract of a $CW$-complex, then it is a cofibration of the underlying topological spaces. 
\end{proposition}
\begin{proof}(sketch) 
We will use the same notations as in lemma \ref{cofmon}. Let $X\rightarrow Y$ a cofibration of CW-complexes and let $\F(X)\rightarrow G$ be any homomorphism of topological groups, then we have an induced homomorphism of topological groups $G\rightarrow \F(Y)\ast_{\F(X)}G$ factors as $G\rightarrow \M(Y)\ast_{\M(X)}G\rightarrow \F(Y)\ast_{\F(X)}G$ in the category of topological monoids. The first map is a cofibration of underlying space \ref{cofmon} and the second map $\M(Y)\ast_{\M(X)}G\rightarrow \F(Y)\ast_{\F(X)}G$ is obtained by attaching cells of underlying spaces. Therefore $G\rightarrow \F(Y)\ast_{\F(X)}G$ is a cofibration of underlying spaces. Using the same procedure as in the proof of \ref{cofmon}, we conclude that if $G\rightarrow H$ is a cofibration of topological groups, then it is a cofibration of the underlying topological spaces.

\end{proof}

\begin{proposition}\label{cw-complex}

Cofibrant topological groups and cofibrant topological abelian groups are Hausdorff.

\end{proposition}

\begin{proof}
We consider the Quillen adjunction between simplicial groups and topological groups
\begin{equation}
\xymatrix{\mathsf{sGr}\ar@<2pt>[r]^-{|-|} & \mathsf{TGr} \ar@<2pt>[l]^-{\mathbf{Sing}} }
 \end{equation} 

and similarly the Quillen adjunction between simplicial abelian groups and topological abelian groups

\begin{equation}
\xymatrix{\mathsf{sAb}\ar@<2pt>[r]^-{|-|} & \mathsf{TAbGr} \ar@<2pt>[l]^-{\mathbf{Sing}} }
 \end{equation}

We use two well known facts:
\begin{enumerate}
\item If $X$ is a topological space, then the counit map $|\mathbf{Sing}(X)|\rightarrow X $ is a fibration. 
 \item If $X_{\bullet}\rightarrow Y_{\bullet}$ is a fibration of simplicial sets, then the geometric realization 
$|X_{\bullet}|\rightarrow |Y_{\bullet}|$ is a fibration. 
\end{enumerate}
Now, suppose that $Z$ is a cofibrant topological group, then in the category of simplicial groups there is a cofibrant simplicial group $G_{\bullet}$ and a trivial fibration $G_{\bullet}\rightarrow \mathbf{Sing}Z$ of simplicial groups, in particular $|G_{\bullet}|\rightarrow |\mathbf{Sing}Z|$ is a trivial fibration of topological groups where $|G_{\bullet}|$ is a cofibrant topological group. Therefore the composition $|G_{\bullet}|\rightarrow |\mathbf{Sing}Z|\rightarrow Z$ is a trivial fibration between cofibrant topological groups. Applying the standard property of left lifting property, it follows that we have a section $s: Z\rightarrow |G_{\bullet}|$ in the category of topological groups. In other words, $Z$ is a retract of $|G_{\bullet}|$. Hence $Z$ is Hausdorff. 

The same argument applies in the abelian case. And we have that any cofibrant topological abelian group is a retract of some $|G_{\bullet}|$ were $G_{\bullet}$ is a cofibrant simplicial abelian group. In particular any cofibrant topological abelian group is Hausdorff. 
\end{proof}

\section{The functors $\E$ and $\Lc$}\label{section2}

\begin{notation}
We fix our notations. 
\begin{enumerate}
\item $R$ is the abelian Lie group of real numbers and $R^{\delta}$ is the additive topological group or real number with discrete topology.
\item The Lie group (with standard topology) $SO(2)$ is denoted by $T$ and the topological group $T^{\delta}$ is the group $T$ with discrete topology.
\item For any topological group $H$, the slice category $H/\Gr$ is defined as follows:
\begin{itemize} 
\item Objects are $H\rightarrow D$ (morphisms of topological groups) 
\item Morphisms between $f:H\rightarrow V$ and $g:H\rightarrow S$ is a morphism of topological groups $h:V\rightarrow S$ such that $h\circ f= g $.
\end{itemize}
\item For any topological group $H$, the category $H-\Top$ is defined as follows:
\begin{itemize} 
\item Objects are topological spaces with a continuous action of $H$.
\item Morphisms are morphisms of topological spaces which commute with the group action of $H$. 
\end{itemize}
\end{enumerate}
\end{notation}
The natural morphism of topological groups $R^{\delta}\rightarrow R$ (identity in the category of groups) induces two kinds of adjunctions that we will consider later:
\begin{equation}
\xymatrix{R^{\delta}/\Gr\ar@<2pt>[r]^-{\Lc} & R/\Gr \ar@<2pt>[l]^-{\U} }
\label{adj2} 
\end{equation}
and
\begin{equation}
\xymatrix{R^{\delta}-\Top\ar@<2pt>[r]^-{\E} & R-\Top \ar@<2pt>[l]^-{\U} }
\label{adj3} 
\end{equation}
The functor $\U$ is the forgetful functor, and for every object $R^{\delta}\rightarrow G$ in $R^{\delta}/\Gr$ we define $\Lc(G)$ by $R\rightarrow R\ast_{R^{\delta}} G$, where $R\ast_{R^{\delta}} G$ is the pushout of the diagram $G\leftarrow R^{\delta}\rightarrow R$ in the category $\Gr$. The functor $\E $ is given by $\E(Y)=T\times_{R^{\delta}} Y.$

\begin{remark}
The forgetful functor $\U$ in \ref{adj2} and \ref{adj3} commutes with limits and colimits. 
\end{remark}
\begin{proposition}\label{idempotent}
The functor $\U\Lc:R^{\delta}/\Gr\longrightarrow R^{\delta}/\Gr$ is idempotent (we will denote $\U\Lc$ again by $\Lc$).
\end{proposition}
\begin{proof}
Let $R^{\delta}\rightarrow G$ be an object of $R^{\delta}/\Gr$, applying the functor $\Lc$ we obtain $R\rightarrow R\ast_{R^{\delta}}G$ which is an object of 
$R/\Gr$, thus an object of $R^{\delta}/\Gr$ via the forgetful functor i.e., $R^{\delta}\rightarrow T\rightarrow R\ast_{R^{\delta}}G$ is an object of  $R^{\delta}/\Gr$, applying again the functor $\Lc$ we obtain $R\rightarrow R\ast_{R^{\delta}} R\ast_{R^{\delta}}G$ but $R\ast_{R^{\delta}} R= R$ since by definition of the pushout in $\Gr$ we have that $R\rightarrow R\ast_{R^{\delta}} R \rightarrow R$ is the identity map in the category of topological groups. Hence $ R\ast_{R^{\delta}}R\ast_{R^{\delta}}G= R\ast_{R^{\delta}}G$. Therefore $\U\Lc$ is idempotent.  

\end{proof}
\begin{notation}
In what follows, we denote an object $R^{\delta}\rightarrow G$ of $R^{\delta}/\Gr$ simply by $G$.   
\end{notation}

\begin{proposition}\label{pushout1}
Let $R^{\delta}\rightarrow G$ and $R^{\delta}\rightarrow H$ be objects of the category $R^{\delta}/\Gr$, the following square in the category of topological groups 
$$ \xymatrix{G\ar[r]\ar[d] & \Lc(G)\ar[d]\\
H\ar[r]& \Lc(H)}$$
 is a pushout. 
\end{proposition}
\begin{proof}
By definition the pushout is given by $\Lc(G)\ast_{G}H$ or equivalently by $R\ast_{R^{\delta}} G\ast_{G}H$ which is isomorphic to $R\ast_{R^{\delta}}H=\Lc(H)$.
\end{proof}

\begin{remark}
The category of $R^{\delta}/\Gr $ is a cofibrantly generated left proper topological model category where equivalences, cofibrations and fibrations are the equivalences, cofibrations and fibrations of the underlying model category $\Gr$ (cf. \ref{model}). The proof can be found in \cite[Theorem 2.8]{Hirs05}.
\end{remark}
\begin{theorem}\label{cartesian}
Suppose that $H$ is a topological group and $R^{\delta}$ a central subgroup in $H$ then $\E(H)=\Lc(H)$. 
\end{theorem}
\begin{proof}
First of all we prove that $\E(H)$ is a topological group. By definition $\E(H)=R\times_{R^{\delta}}H$ where $R^{\delta}$ acts on $R$ by right translation and on $H$ by left translation. Let consider the continuous map 
$$(R\times H)\times (R\times H)\rightarrow R\times H$$
$$(t,h), (t^{'},h^{'})\mapsto (tt^{'},hh^{'}) $$
Since $R^{\delta}$ is central in $H$, then for any $a,b\in R^{\delta}$ we have that  
$$(ta,a^{-1}h), (t^{'}b,b^{-1}h^{'})\mapsto (tt^{'}ab,(ab)^{-1}hh^{'}) $$
in particular we have a continuous induced map
$$\E(H)\times \E(H)\rightarrow \E(H)$$ 
endowing $\E(H)$ with the structure of topological monoid. In the same way the map 
$$inv:R\times H\rightarrow R\times H$$
$$(t,h)\mapsto (t^{-1},h^{-1})$$ 
factors to  
$$inv:\E(H)\rightarrow \E(H)$$
since $R^{\delta}$ is central in $H$. In particular the evident continuous bijective map $H\rightarrow \E(H)$ is a map of topological groups. Let consider the homomorphism of groups $H\rightarrow \Lc(H)$, in particular it is a homomorphism of $R^{\delta}$-spaces, hence it factors as $H\rightarrow\E(H)\rightarrow \Lc(H)$. By the universal property of the pushout, we conclude that $\E(H)=\Lc(H)$ in the category of topological groups. 
\end{proof}
\subsection{Topological description of $\Lc(G)$ }
Let $H$ be a topological group and let $R^{\delta}$ a topological subgroup of $H$ (for simplicity). In this paragraph, we study the following pushout 
 $$\xymatrix{R^{\delta}\ar[r]^{i}\ar[d] & H\ar[d]\\
R\ar[r] & \Lc(H) \cocartesien
} $$ 
in the category $\Gr$, more precisely we study the topology of $\Lc(H)$.
\begin{definition}\label{defpush}
Let $H$ be a topological group and $R^{\delta}$ a topological subgroup of $H$ with multiplication $\mu$, we define $F^{(i)}(H)$ as follows 
\begin{enumerate}
\item $F^{(0)}(H)=H$
\item $F^{(1)}(H)$ is given by the pushout diagram in the category of topological spaces 
$$\xymatrix{ H\times R^{\delta}\times H\ar[r]^-{\mu}\ar[d] & F^{(0)}(H)\ar[d]\\
H\times R\times H\ar[r]^-{\mu} & F^{(1)}(H) \cocartesien
} $$
\item  $F^{(2)}(H)$ is given by the pushout diagram in the category of topological spaces 
$$\xymatrix{ & H\ar[d] &\\
H\times R^{\delta}\times H\times R^{\delta}\times H\ar[r]^-{\mu}\ar[d]\ar[ru]^-{\mu} & F^{(1)}(H)\ar[d]\\
H\times R\times H\times R\times H\ar[r]^-{\mu} & F^{(2)}(H) \cocartesien
} $$
\item $\dots$
\end{enumerate}
\end{definition}
\begin{lemma}\label{push1}
Under the same notation as in definition \ref{defpush}, we have that the colimit in the category of topological spaces 
$$F^{(\infty)}(H):=\colim[F^{(0)}(H)\rightarrow F^{(1)}(H)\rightarrow F^{(2)}(H)\rightarrow \dots]=\Lc(H) $$  
\end{lemma}
\begin{proof}
First of all, all maps in the diagrams that are not labelled are identity maps of underlying sets. We recall that in the category of compactly generated spaces the cartesian product commutes with colimits. We notice that the natural map of topological groups $H\rightarrow \Lc(H)$ factors through $F^{(\infty)}(H)$. It can be checked step by step, the pushout diagram in the category of topological spaces 
$$\xymatrix{ H\times R^{\delta}\times H\ar[r]^-{\mu}\ar[d] & F^{(0)}(H)\ar[d]\ar[rdd]\\
H\times R\times H\ar[r]^-{\mu}\ar[d] & F^{(1)}(H)\cocartesien \ar@{.>}[rd]^-{\exists !} \\
 \Lc(H)\times R\times \Lc(H) \ar[rr]^-{\Lc(\mu)}& & \Lc(H)} $$ 
 implies the existence of a unique continuous map $F^{(1)}(H)\rightarrow \Lc(H)$. By induction we have a pushout diagram 
 $$\xymatrix{ H\times R^{\delta}\times H\dots \times R^{\delta}\times H\ar[r]^-{\mu}\ar[d] & F^{(n)}(H)\ar[d]\ar[rdd]\\
H\times R\times H\dots \times R\times H\ar[r]^-{\mu}\ar[d] & F^{(n+1)}(H)\cocartesien \ar@{.>}[rd]^-{\exists !} \\
 \Lc(H)\times R\times \Lc(H)\dots \times R\times \Lc(H) \ar[rr]^-{\Lc(\mu)}& & \Lc(H)} $$ 
It follows that $H\rightarrow \Lc(H)$ factors through $F^{(\infty)}(H)$ in the category of spaces and in the category of abstract groups since as groups $H, \Lc(H)$ and $F^{(\infty)}(H)$ are the same. 
By construction the multiplication $\mu:  F^{(n)}(H)\times  F^{(m)}(H)\rightarrow  F^{(n+m)}(H)$ is continuous, therefore 
the induced multiplication 
$\mu:\colim  F^{(n)}(H)\times  \colim F^{(m)}(H)\rightarrow  \colim F^{(n+m)}(H)$
is continuous, hence $\mu: F^{(\infty)}(H)\times F^{(\infty)}(H)\rightarrow F^{(\infty)}(H)$ is continuous. For instance we can illustrate the continuity in the case where $n=m=1$, by definition the map $(H\times R\times H) \times (H\times R\times H)\rightarrow F^{(2)}(H)$ defined as the composition 
$$(h_{1},t_{1},h_{2},h_{3},t_{2},h_{4})\mapsto (h_{1},t_{1},\mu(h_{2},h_{3}),t_{2},h_{4})\mapsto \mu(h_{1},t_{1},\mu(h_{2},h_{3}),t_{2},h_{4})$$ is continuous
and for the same reason the map $$(H\times R\times H) \times [(H\times R\times H)\sqcup F^{(0)}(H)]\rightarrow F^{(2)}(H)$$ is continuous 
hence the induced map $$(H\times R\times H) \times [(H\times R\times H)\sqcup_{H\times R^{\delta}\times H} F^{(0)}(H)]\rightarrow F^{(2)}(H)$$ 
is continuous therefore $(H\times R\times H) \times F^{(1)}(H)\rightarrow F^{(2)}(H)$ is continuous, and by the same argument $\mu: F^{(1)}(H) \times F^{(1)}(H)\rightarrow F^{(2)}(H)$ is continuous.
Following the same method we obtain that $\mu:F^{(k)}(H)\times F^{(m)}(H)\times F^{(n)}(H)\rightarrow F^{(k+m+n)}(H)$ is an associative continuous multiplication and that the inverse map $(h\mapsto h^{-1})$ given by 
$$ inv:F^{(\infty)}(H)\rightarrow F^{(\infty)}(H)$$ is continuous since at each stage $inv:F^{(n)}(H)\rightarrow F^{(n)}(H)$ is continuous . Hence $F^{(\infty)}(H)$ is a topological group. We have proved that the factorisation 
$$H\rightarrow F^{(\infty)}(H)\rightarrow  \Lc(H)$$ in the category of topological spaces is a factorisation in the category of topological groups. By the universality of the pushout, we conclude that $F^{(\infty)}(H)=\Lc(H)$. 

\end{proof}

\begin{remark}
It can be proven that $F^{(n)}(H)$ is canonically homeomorphic to the space defined as the pushout 
$$\xymatrix{ H\times R^{\delta}\times H\times \dots R^{\delta}\times H \ar[rrr]^{\mu} \ar[d] &&& H\ar[d]\\
H\times R\times H\times \dots R\times H \ar[rrr]^{\mu} &&& F^{(n)}(H)\cocartesien
 } $$
in the category of topological spaces. 
\end{remark}
\subsection{Alternative computation of $\Lc(G)$ }\label{section4bis}

Let $R^{\delta}$ be a topological subgroup of $H$, in this section we give an other useful alternative to compute the pushout digram $R\ast_{R^{\delta}} H$ for any  object $H\in R^{\delta}/\Gr$. We recall that $\E(H):=R\times_{R^{\delta}} H$. 

\begin{definition}
Let $H\in R^{\delta}/\Gr$, a topological group with multiplication $\mu$, and a homomorphism of topological groups $R^{\delta}\rightarrow H$ (for simplicity the homomorphism is an embedding). We define $\E^{n}(H)$ as the pushout in the category of topological spaces:

$$\xymatrix{H^{\times^{n}}\ar[r]^{\mu}\ar[d]^{id} & H\ar[d]_{id}\\
\E (H)^{\times^{n} }\ar[r] & \E^{n}(H)} $$
where $H^{\times^{n}}$ is the $n$-times cartesian product of $H$. 
\end{definition}
 
From the definition it is clear that there is a continuous "identity" map $\E^{n}(H)\rightarrow \E^{n+1}(H)$.
\begin{proposition}\label{alt}
The filtered colimit $$\E^{\infty}(H):= colim [\E(H)\rightarrow \dots \rightarrow \E^{n}(H)\rightarrow \E^{n+1}(H)\rightarrow\dots ]$$ in the category of topological spaces is exactly $\Lc(H)$ (as space).   
\end{proposition} 
\section{Left properness}\label{section3}

In this section we give a proof that the model category  $\Gr$ defined in \ref{model} is left proper. 
\begin{definition}
A model category is left proper if the pushout of any weak equivalence along any cofibration is a weak equivalence.  
\end{definition} 

In order to prove the theorem \ref{proper} we need to introduce some known results. The category of topological monoids $\mathsf{TMon}$ \ref{grmon} is in fact a cofibrantly generated model category. The adjunction \ref{adj1} can be decomposed in a pair of Quillen adjunctions:

\begin{equation}
\xymatrix{\Top\ar@<2pt>[r]^-{\M} & \mathsf{TMon} \ar@<2pt>[l]^-{\U}\ar@<2pt>[r]^-{\C} & \Gr \ar@<2pt>[l]^-{\U}}
\label{adj4} 
\end{equation}
Where the composition $ \C\circ \M:=\C\M$ is the functor $\F$ defined in the adjunction \ref{adj1}. We notice that the functor $\C$ is sometimes called the (naive) completion functor. The generating cofibrations in the cofibrantly generated model category of topological monoids $\mathsf{TMon}$ are given by $\M(\partial \Delta^{n})\rightarrow \M(\Delta^{n})$ with $n\in\mathbf{N}$.
\begin{theorem}\cite[3.1. Theorem (b'), example 1.19]{batanin2017homotopy}\label{bb}
The model category of topological monoids $\mathsf{TMon}$ is left proper. 
\end{theorem}
\begin{theorem}\label{completionlemma}
Suppose that $A$ is a topological group and suppose that $X\rightarrow Y$ is a cofibration in the category of topological sapces in particular $\M(X)\rightarrow \M(Y)$ (resp. $\F(X)\rightarrow \F(Y)$) is a cofibration in $\mathsf{TMon}$  (resp. in $\Gr$). Then the natural map of topological monoids 
$$A\ast_{\M(X)}\M(Y)\rightarrow A\ast_{\M(X)}\F(Y)= A\ast_{\F(X)}\F(Y)$$ 
induces a weak equivalence of topological spaces $$\mathbf{B}[A\ast_{\M(X)}\M(Y)]\rightarrow \mathbf{B}[ A\ast_{\F(X)}\F(Y)].$$ 
\end{theorem}
\begin{proof}

Let $A$ be a topological group. We consider the slice categories $A/\Gr$ and $A/\mathsf{TMon}$ these two categories are naturally model categories. Since $A$ is a topological group $A=\C(A)$, we have a Quillen adjunction 
$$\xymatrix{A/\mathsf{TMon}\ar@<2pt>[r]^-{\C} & A/\Gr \ar@<2pt>[l]^-{\U} }$$
Suppose that $A\rightarrow D$ is a cofibrant object in  $A/\mathsf{TMon}$ and  $A\rightarrow H$ any object in $A/\Gr$,  then by the previous adjunction $\Map_{A/\mathsf{TMon}}(D,H)=\Map_{A/\Gr}(\C(D),H)$. By theorem \ref{classifying},
we have that  $$\Map_{A/\mathsf{TMon}}(D,H)\simeq \Map_{\B A/\Top_{\ast}}(\B D,\B H)$$ and 
$$\Map_{A/\Gr}(\C(D),H)\simeq \Map_{\mathbf{B}A/\Top_{\ast}}(\mathbf{B}\C(D),\mathbf{B}H)$$
therefore, the map $D\rightarrow \C(D)$ induces an equivalence of topological spaces $$ \Map_{\mathbf{B}A/\Top_{\ast}}(\B\C(D),\B H)\rightarrow \Map_{\B A/\Top_{\ast}}(\B D,\B H)$$ for any object $H\in A/\Gr$. In particular $\mathbf{B}D\rightarrow \B\C(D)$ is an equivalence of topological spaces. Noticing that $\C(A\ast_{\M(X)}\M(Y))=A\ast_{\F(X)}\F(Y)$ we prove our theorem. 
\end{proof}

\begin{lemma}\label{paretape}
Suppose we have a pushout diagram in $\Gr$ :
\begin{equation}\label{diaggroup}
\xymatrix{A\ar@{>->}[d]\ar[r]^{\sim} & C\ar[d]^-{g}\\
A\ast_{\F(X)}\F(Y) \ar[r]& \cocartesien C \ast_{\F(X)}\F(Y)  }
\end{equation}
where the left vertical map is a cofibration and the upper horizontal map is a weak equivalence moreover we assume that $\F(X)\rightarrow \F(Y)$ is a generating cofibration. Then the lower horizontal map $$A\ast_{\F(X)}\F(Y)\rightarrow  C \ast_{\F(X)}\F(Y) $$ is a weak equivalence of topological groups.  

\end{lemma}
\begin{proof}
We consider the following pushout diagram in the category $\mathsf{TMon}$ where $\C(m)=g$ and $A\rightarrow A\ast_{\M(X)}\M(Y)$ is a cofibration in the category $\mathsf{TMon}$:
\begin{equation}
\xymatrix{A\ar@{>->}[d]\ar[r]^{\sim} & C\ar[d]^-{m}\\
A\ast_{\M(X)}\M(Y) \ar[r]& \cocartesien C \ast_{\M(X)}\M(Y)  }
\end{equation}
Then combining this diagram with \ref{diaggroup} we obtain a commutative diagram in $\mathsf{TMon}$: 
\begin{equation}
\xymatrix{A\ar[r]^{id}\ar[d]& A\ar[d]\ar[r]^{\sim} & C\ar[d]^-{g}\ar@/^3pc/[dd]^-{m}\\
A\ast_{\M(X)}\M(Y) \ar[r]^{\sim_{\B}} \ar[d]^{id}& A\ast_{\F(X)}\F(Y) \ar[r]& C \ast_{\F(X)}\F(Y)\cocartesien\\  
A\ast_{\M(X)}\M(Y) \ar[rr]^{\sim} & & C \ast_{\M(X)}\M(Y)\ar[u]^{\sim_{\B}} \cocartesien}
\end{equation}
By theorem \ref{bb}, the map $ A\ast_{\M(X)}\M(Y)\rightarrow C\ast_{\M(X)}\M(Y) $ is a weak equivalence of topological monoids and by theorem \ref{completionlemma} we have that 
$\mathbf{B}[A\ast_{\M(X)}\M(Y)]\rightarrow \mathbf{B}[A\ast_{\F(X)}\F(Y)]$ and 
$\mathbf{B}[C\ast_{\M(X)}\M(Y)]\rightarrow \mathbf{B}[C\ast_{\F(X)}\F(Y)]$ are weak equivalences of topological spaces. We conclude that the map 
$ \mathbf{B}[A\ast_{\F(X)}\F(Y)]\rightarrow \mathbf{B}[C\ast_{\F(X)}\F(Y)] $ is a weak equivalence hence $ A\ast_{\F(X)}\F(Y)\rightarrow C\ast_{\F(X)}\F(Y) $
is a weak equivalence of topological groups. 
\end{proof}
\begin{theorem}\label{proper}
The category of topological groups is left proper.
\end{theorem}
\begin{proof}
Suppose we have a pushout diagram in $\Gr$ 
\begin{equation}\label{diag1}
\xymatrix{A\ar[d]\ar[r]^{\sim} & C\ar[d]\\
B \ar[r]& D  \cocartesien }
\end{equation}
 By Lemma \ref{paretape}, if the cofibration $A\rightarrow B$ is of the form $A\rightarrow A\ast_{\F(X)}\F(Y)$ then the map $B\rightarrow D$ is a weak equivalence. If the cofibration $A\rightarrow B$ is cell-cofibration \ref{cellcof} then by the same argument \ref{paretape} we have that $B\rightarrow D$ is a weak equivalence (since cell-cofibration are defined by transfinite composition ). But any cofibration in $\Gr$ is retract of such of a cell-cofibration, therefore $B\rightarrow D$ is a weak equivalence for any cofibration $A\rightarrow B$. Hence $\Gr$ is a left proper model category. 
\end{proof}

\section{The case of central extensions}\label{section4}

This section is independent form the rest, and can be taken as a warm-up for the results described in the section \ref{section5}.  
\begin{theorem}\label{central}
Suppose that $H$ is topological group and $R^{\delta}$ is a topological subgroup of $H$ such that
\begin{enumerate}
\item $R^{\delta}$ is central in $H$,
\item as a topological space $H$ is path connected and Hausdorff,
\end{enumerate} 
Then the induced map $$\mathrm{H}_{\ast}(H;\mathbf{F}_{p})\rightarrow \mathrm{H}_{\ast}(\Lc(H);\mathbf{F}_{p}) $$ 
is an isomorphism. 
\end{theorem}
\begin{proof}
Since $R^{\delta}$ is central in $H$, by \ref{cartesian} we have that $\E(H)=\Lc(H)$. By hypothesis $H$ is connected and Hausdorff and the action (by translation) of $R^{\delta}$ on $H$ is proper and free, therefore we have a Cartan Leray spectral sequence \cite[Theorem $8^{bis}9$]{mccleary2001user} which the second page converges to the homology of $\E(H)$ given by:
$$E_{n,m}^{2}=\mathrm{H}_{n}(R^{\delta};\mathrm{H}_{m}(R\times H;\mathbf{F}_{p}))\Longrightarrow \mathrm{H}_{n+m}(\E(H);\mathbf{F}_{p}). $$
Since the topological group $H$ is path connected, the action (by left multiplication) of $R^{\delta}$ on $H$ is homotopy equivalent to the identity. In particular the action of action of $R^{\delta}$ on $\mathrm{H}_{m}(H;\mathbf{F}_{p})$ is trivial, therefore 
$$E_{n,m}^{2}=\mathrm{H}_{n}(\B R^{\delta};\mathrm{H}_{m}(R\times H;\mathbf{F}_{p}))\Longrightarrow \mathrm{H}_{n+m}(\E(H);\mathbf{F}_{p}),$$ 
or more simply
$$E_{n,m}^{2}=\mathrm{H}_{n}(\B R^{\delta};\mathrm{H}_{m}(H;\mathbf{F}_{p}))\Longrightarrow \mathrm{H}_{n+m}(\E(H);\mathbf{F}_{p}).$$ 
But the homology groups $\tilde{\mathrm{H}}_{\ast}(H;\mathbf{F}_{p})$ are $p$-torsion it follows by \cite[Lemma 1.1]{suslin1983thek}, that the spectral sequence $E_{n,m}^{2}$ collapses (since  the reduced integral homology groups of $\B R^{\delta}$ are rational vector spaces). In particular  
$$ \mathrm{H}_{\ast}(\E(H);\mathbf{F}_{p}) = \mathrm{H}_{\ast}(H;\mathbf{F}_{p})$$

Hence the natural map $H\rightarrow \E(H)=\Lc(H)$ induces an isomorphism in homology with coefficients in $\mathbf{F}_{p}$. 
\end{proof}

\begin{remark}
The reason we have mentioned theorem \ref{central} is that, a priori, it can not be recovered from theorem \ref{crucial+}.
\end{remark}



\section{Comparison map $G\rightarrow \Lc (G)$ }\label{section5}


\begin{lemma}\label{bb-pushout}
Let $X_{\bullet}\rightarrow Y_{\bullet}$ and $X_{\bullet}\rightarrow Z_{\bullet}$ be injective homomorphisms of simplicial groups. Then le following pushout square in the category of simplicial groups 
\begin{equation}\label{eqation12}
\xymatrix{ X_{\bullet}  \ar[r] \ar[d] & Y_{\bullet}\ar[d]\\
X_{\bullet}  \ar[r]  & D_{\bullet}\cocartesien
} 
\end{equation}
is a homotopy pushout
\end{lemma}
\begin{proof}
 First, we factor the map $X_{\bullet}\rightarrow Z_{\bullet}$ as cofibration followed by a trivial fibration $X_{\bullet}\rightarrow Z_{\bullet}^{'}\rightarrow Z_{\bullet}$ in the model category of simplicial groups. Let consider the double pushout diagram in the category of simplicial groups given by 
 
 $$\xymatrix{ X_{\bullet}  \ar[r] \ar@{>->}[d] & Y_{\bullet}\ar[d]\\
 Z_{\bullet}^{'}\ar@{->>}[d]^{\sim} \ar[r] & D_{\bullet}^{'}\ar[d]\cocartesien\\
Z_{\bullet}  \ar[r]  & D_{\bullet}\cocartesien
} $$
Our goal is to prove that the homomorphism $D_{\bullet}^{'}\rightarrow D_{\bullet}$ is a weak equivalence of simplicial groups (i.e. weak homotopy equivalence of underlying simplicial sets). Applying the simplicial classifying functor, we have that the two diagrams in the category of simplicial sets 
$$\xymatrix{ \B_{\bullet} X_{\bullet}  \ar[r] \ar[d] & \B_{\bullet} Y_{\bullet}\ar[d]\\
 \B_{\bullet} Z_{\bullet}  \ar[r]  & \B_{\bullet} D_{\bullet}
} $$
 and 
 $$\xymatrix{ \B_{\bullet}X_{\bullet}  \ar[r] \ar[d] & \B_{\bullet}Y_{\bullet}\ar[d]\\
 \B_{\bullet}Z_{\bullet}^{'} \ar[r] & \B_{\bullet}D_{\bullet}^{'}
} $$
are homotopy pushout squares, it follows from \cite[Theorem 4.0]{fiedorowicz1984classifying} applied degree by degree and taking the diagonal functor from bisimplicial sets to simplicial sets (cf. \cite{seb2008} and remark \ref{remarkclass}). Therefore, the natural map of simplicial sets $\B_{\bullet} D_{\bullet}^{'}\rightarrow \B_{\bullet} D_{\bullet}$ is a weak homotopy equivalence. Thus the homomorphism $D_{\bullet}^{'}\rightarrow D_{\bullet}$ is a weak equivalence of simplicial groups. But the model category of simplicial groups is left proper since the category of topological groups is left proper. Thus the diagram \ref{eqation12} is a homotopy pushout square.
\end{proof}
\begin{theorem}\label{B-pushout}
Suppose that $G\rightarrow X$ is a cofibration of topological groups and $G\rightarrow H$ is an injective homomorphism of topological groups, then $$\xymatrix{ \B G\ar[r]\ar[d] & \B X\ar[d]\\
\B H\ar[r] & \B [H\ast_{G}X ]
} $$
 is a homotopy pushout of topological spaces. In particular if $X$ is contractible as a space, $G=R^{\delta}$ and $H=R$, then the topological group $R\ast_{R^{\delta}}X$ is weakly homotopy equivalent, as a space, to $\Omega\Sigma\B R^{\delta}$.

\end{theorem}
\begin{proof}
Let consider the pushout diagram in the category of topological groups 
$$\xymatrix{ G\ar@{>->}[r]\ar@{(->}[d] & X\ar@{(->}[d]\\
H\ar@{>->}[r] & H\ast_{G}X \cocartesien
} $$
Since the model category of topological groups is left proper, the previous pushout diagram is a homotopy pushout square. On the other hand the induced maps of simplicial groups $\Sing G\rightarrow  \Sing X$ and $\Sing G\rightarrow  \Sing H$
are monomorphisms, therefore the pushout diagram in the category of simplicial groups 
$$\xymatrix{ \Sing G\ar[r]\ar[d] & \Sing X\ar[d]\\
\Sing H\ar[r] & \Sing H\ast_{ \Sing G}\Sing X \cocartesien
} $$ is a homotopy pushout  square of simplicial groups (cf \ref{bb-pushout}), in particular the pushout diagram in the category of topological groups 
$$\xymatrix{ |\Sing G|\ar[r]\ar[d] & |\Sing X|\ar[d]\\
|\Sing H|\ar[r] & |\Sing H\ast_{ \Sing G}\Sing X|\cocartesien }$$

is a homotopy pushout in the category of topological groups, it follows that the natural map of topological groups $|\Sing H\ast_{ \Sing G}\Sing X|\rightarrow H\ast_{G}X$ is a weak homotopy equivalence, on the other hand applying the simplicial classifying functor $\B_{\bullet}$ we obtain that 
$$\xymatrix{ |\B_{\bullet}\Sing G|\ar[r]\ar[d] & |\B_{\bullet}\Sing X|\ar[d]\\
|\B_{\bullet}\Sing H|\ar[r] & |\B_{\bullet}[\Sing H\ast_{ \Sing G}\Sing X]| }$$
is a homotopy pushout in the category of topological spaces (we apply \cite[theorem 4.0.]{fiedorowicz1984classifying} degree by degree and \cite{seb2008}), in particular the following commutative square in the category of topological spaces   
$$\xymatrix{ \B G\ar[r]\ar[d] & \B X\ar[d]\\
\B H\ar[r] & \B [H\ast_{G}X ]
} $$
is a homotopy pushout in the category of topological spaces ($|\B_{\bullet}\Sing(-) |:= \B(-)$). As a consequence, if $G=R^{\delta}$ and $H=R$ then $\B [R\ast_{R^{\delta}}X ]$ is homotopy equivalent to $\Sigma \B R^{\delta}$. We conclude that 
$R\ast_{R^{\delta}}X$ is homotopy equivalent to $\Omega \Sigma \B R^{\delta}$  as topological space.

\end{proof}
\subsection{Precofibrant $R^{\delta}$-topological groups}
\begin{definition}\label{qcofibrant}
Let $X$ be an object of $R^{\delta}/\Gr$, we say that $X$ is \textbf{precofibrant} $R^{\delta}$-topological groups if $\Lc X= R\ast_{R^{\delta}}X$ is a cofibrant object in $R/\Gr$.
\end{definition}

\begin{proposition}
Let $R^{\delta}\rightarrow H$ be a precofibrant object in $R^{\delta}/\Gr$, then $1\ast_{R^{\delta}} H$ is a cofibrant object in $\Gr$. Moreover, $X$ is Hausdorff, since $R\rightarrow\Lc X$ is cofibrantion and hence $\Lc X$ is Hausdorff.  
\end{proposition}
\begin{proof}
Since $R^{\delta}\rightarrow H$ be a precofibrant, then by definition $R\rightarrow \Lc H$ is a cofibrant object in $R/\Gr$ or equivalently, $R\rightarrow \Lc H$ is a cofibration in $\Gr$. By definition $\Lc H= R\ast_{R^{\delta}}H$, thus $1\ast_{R^{\delta}} H =1\ast_{R} \Lc H$, but $1\ast_{R} \Lc H$ is cofibrant in $\Gr$ (pushout of a cofibration is a cofibration). Since $R\rightarrow \Lc H$ is a cofibration in $\Gr$, it follows that $\Lc H$ is Hausdorff (cf. \ref{coftgroup}, \ref{cw-complex}), therefore $H$ is Hausdorff.
\end{proof}

\begin{definition}\label{ratdef}
A (path) connected topological group $X$ is \textbf{rational} if all homotopy groups $\pi_{\ast}(X)$ are rational vector spaces.
An object $R^{\delta}\rightarrow X$ of $R^{\delta}/\Gr$ is \textbf{rational} if  $X$ is path connected and all homotopy groups $\pi_{\ast}(X)$ are rational vector spaces.
\end{definition}
\begin{proposition}\label{rat1}
A path connected topological group $X$ is rational if and only if all reduced homology groups $\tilde{\mathrm{H}}_{\ast}(X;\mathbf{Z})$ are rational vector spaces \cite[definition 1.1]{hess2007rational}.
\end{proposition}

\begin{proposition}\label{abelianizationfun}
The abelianization functor $$\mathbf{Ab}: \Gr\rightarrow \mathsf{AbTGr} $$
from the category of topological groups to the category of abelian topological groups is the left adjoint to the forgetful functor (in particular it commutes with colimits).
\end{proposition}

\begin{proposition}\label{rationality0}
If $Z$ is a cofibrant topological group then $$\pi_{\ast}\mathbf{Ab}(Z)\cong \mathrm{H}_{\ast+1}(\B Z;\mathbf{Z}).$$ 
\end{proposition}
\begin{proof}

Let $\mathbf{Sing}Z:=Z_{\bullet}$ the simplicial group associated to $Z$. Since the category of simplicial groups $\mathsf{sGr}$ is  a model category, there exists a trivial fibration of simplicial groups  $G_{\bullet}\rightarrow Z_{\bullet}$ such that $G_{\bullet}$ is a cofibrant simplicial group. Therefore the composition $|G_{\bullet}|\rightarrow |Z_{\bullet}|$ with the counit natural transformation $|Z_{\bullet}|\rightarrow Z$ is a weak equivalence. By construction, it follows that
\begin{itemize} 
\item we have weak homotopy equivalences between  $\mathbf{B}|G_{\bullet}|\sim |\mathbf{B}G_{\bullet}|\sim \mathbf{B}|Z_{\bullet}|\sim \mathbf{B} Z$. 
\item $|\mathbf{Ab}(G_{\bullet})|$ is isomorphic to $\mathbf{Ab}(|G_{\bullet}|)$ as abelian topological groups.
\item the induced map $|\mathbf{Ab}(G_{\bullet})|=\mathbf{Ab}(|G_{\bullet}|)\rightarrow \mathbf{Ab}(Z)$ is a weak homotopy equivalence since $\mathbf{Ab}$ is a left Quillen functor and $|G_{\bullet}|$ and $Z$ are cofibrant topological groups by hypothesis. 
\end{itemize}
By \cite[example 4.26]{goerss2007model}, we have that
$$\pi_{\ast}(\mathbf{Ab}(G_{\bullet}))\cong \mathrm{H}_{\ast+1}(\mathbf{B}G_{\bullet});\mathbf{Z} ) $$
which implies  
$$\pi_{\ast}(\mathbf{Ab}(Z))\cong \mathrm{H}_{\ast+1}(\mathbf{B}Z);\mathbf{Z} ) $$ 
\end{proof}
\begin{proposition}
If $Z$ is a topological (path) connected topological group such that for any prime $p$ the $p$-power maps $\mathbf{p}:Z\rightarrow Z$, $z\mapsto z^{p}$  is a weak homotopy equivalence of underlying spaces, then $Z$ is rational. Moreover,
 $$\pi_{\ast}\mathbf{p}:\pi_{\ast}Z\rightarrow \pi_{\ast}Z$$ $$ [\alpha]\mapsto p[\alpha] .$$ 
\end{proposition}
\begin{proof}
Suppose we have two loops $\gamma, \sigma: [0,1]\rightarrow Z$, then up to homotopy the concatenation of loops $[\gamma\star \sigma]$ represent the same class as $[\gamma \cdot \sigma]$ where $\cdot$ is the pointwise multiplication ($Z$ is a topological group). For  $n$-th homotopy groups, we use the same argument for loops in $\Omega^{n-1}Z$. The result follows immediately. If for any prime $p$, the map $\mathbf{p}_{\ast}:\pi_{\ast}Z\rightarrow \pi_{\ast}Z$ is an isomorphism, it implies that the groups $\pi_{\ast}Z$ are rational vector spaces, hence $Z$ is a rational topological group in the sense of \ref{ratdef}.
\end{proof}
In what follows, we will use frequently the following fact: If $\Gamma\subset H$ is a discrete closed subgroup of a Hausdorff topological group $H$, then $H$ is a principal $\Gamma$-bundle \cite[proposition 14.1.12, example 14.1.13]{tom2008algebraic}. 
\begin{lemma}\label{rationality}
Let $X$ be a precofibrant object of $R^{\delta}/\Gr$ such that
\begin{itemize}
\item $X$ is path connected as a topological space.
\item $R^{\delta} \rightarrow X$ and $R^{\delta} \rightarrow \mathbf{Ab}(X)$ are embeddings of topological groups. 
\end{itemize}
then $\Lc (X)$ is rational if and only if $\mathbf{Ab}(X)$ is rational.  

\end{lemma}
\begin{proof}
First, we notice that as a group $X$ is isomorphic to $R\ast \F(Y)$ where $\F(Y)$ is a free group generated by a set $Y$. 
Lets consider a commutative diagram of topological groups
$$\xymatrix{R^{\delta}\ar@{(->}[r]\ar[d] & X \ar[r]^-{quotient}\ar[d] & 1\ast_{R^{\delta}}X \ar@{=}[d]\\
R\ar@{>->}[r] & \Lc(X) \ar[r]^-{quotient} & 1\ast_{R}\Lc(X)
} $$
since $R\rightarrow \Lc(X)$ is a cofibration of topological groups and since that the category of topological groups is left proper it follows that $\Lc(X)\rightarrow 1\ast_{R^{\delta}}X$ is a weak homotopy equivalence of underlying spaces. On another hand, applying the functor $\mathbf{Ab}$ we obtain a short exact sequences of abelian topological groups 
 $$\xymatrix{R^{\delta}\ar[r] & \mathbf{Ab}(X) \ar[r]^-{quotient} & \mathbf{Ab}(X)/R^{\delta}=[1\times\mathbf{Ab}(X)]/ R^{\delta}} $$
Suppose that $\Lc(X)$ is rational it follows that $1\ast_{R^{\delta}}X$ is rational, but $1\ast_{R^{\delta}}X$ is a path connected cofibrant topological group, therefore $\mathbf{Ab}(X)/R^{\delta}=\mathbf{Ab}(1\ast_{R^{\delta}}X)$ is a rational topological group (cf. \ref{rationality0}). Notice that $1\ast_{R^{\delta}}X$ is hausdorff (cf. \ref{cw-complex}) and $\mathbf{Ab}(X)/R^{\delta}$ is a cofibrant topological abelian group, hence it is Hausdorff (cf \ref{cw-complex}). Since $R^{\delta} \rightarrow \mathbf{Ab}(X)$ is an embeddings of topological groups, we have a short exact sequence of Hausdorff topological groups 
$$\xymatrix{R^{\delta}\ar[r] & \mathbf{Ab}(X) \ar[r] & \mathbf{Ab}(X)/R^{\delta} 
} $$
We notice that $R^{\delta}\rightarrow \mathbf{Ab}(X)$ is a closed embedding. It follows that we have a homotopy fiber sequence $\mathbf{Ab}(X)/R^{\delta}\rightarrow \B R^{\delta}\rightarrow \B \mathbf{Ab}(X)$, thus the homotopy groups of $\mathbf{Ab}(X)$ are rational vector spaces. Conversely, if we suppose that $\mathbf{Ab}(X)$ is rational, it implies that $\mathbf{Ab}(X)/R^{\delta}$ is rational, hence $1\ast_{R^{\delta}}X$ is rational (cf. \ref{rationality0}, \ref{rat1}) and finally $\Lc(X)$ is rational.  
\end{proof}

\begin{theorem}\label{contractible}
Suppose that $O$ is a precofibrant object $(\in R^{\delta}/\Gr )$ in the sense of \ref{qcofibrant} such that as a space $O$ is a contractible and $R^{\delta}\rightarrow O$, $R^{\delta}\rightarrow \mathbf{Ab}(O)$ are embeddings, then $\Lc(O)$ is rational. 
\end{theorem}
\begin{proof}
\textbf{Spet 1:}\\
 We factor $Ab: O\rightarrow \mathbf{Ab}(O)$ in the model category $R^{\delta}/\Gr$ as 
$$\xymatrix{O\ar@{>->}[r]^-{\sim} & O^{'}\ar@{->>}[r]^-{Ab^{'}} & \mathbf{Ab}(O)}$$ 

where the first map is a trivial  cofibration and the second map is a fibration. 
It is enough to prove that $\mathbf{Ab}(O)$ is rational (cf. \ref{rationality}). For any prime $p$, we consider the pullback diagram in the category $\Gr$ (hence a pullback in the category of topological spaces)
$$\xymatrix{O_{\mathbf{p}}\ar[r]\ar[d] & \mathbf{Ab}(O)\ar[d]^-{\mathbf{p}}\\
O^{'}\ar[r]^-{Ab^{'}} & \mathbf{Ab}(O)
}$$

We notice that the map $\mathbf{p}(x)=x^{p}$ is an injective (inj) homomorphism of topological groups (since as a group $\mathbf{Ab}(O)$ is a product of $R$ and a transfinite composition of free abelian group). Then the previous pullback diagram is a homotopy pullback (in the category of spaces). Since the map $\mathbf{p}$ is a rational equivalence i.e. induces isomorphism in homotopy groups after tensoring with rational field, it follows that $O_{\mathbf{p}}$ is path connected. On another hand, by universality of the pullback, the map $\mathbf{p}:O_{\mathbf{p}}\rightarrow O_{\mathbf{p}}$ factors trough $O^{'}$. More precisely we have a commutative diagram in the category of topological spaces given by 
\begin{equation}\label{eqpullbak}
\xymatrix{ O^{'}\ar[drr]^-{Ab^{'}}\ar[ddr]_-{\mathbf{p}} \ar@{.>}[rd]& & \\
& O_{\mathbf{p}}\ar@{->>}[r]\ar[d]^-{inj} & \mathbf{Ab}(O)\ar[d]^-{\mathbf{p}}\\
& O^{'}\ar@{->>}[r]^-{Ab^{'}} & \mathbf{Ab}(O)
}
\end{equation}
The map $$O^{'}\rightarrow O_{\mathbf{p}}\subset O^{'}\times \mathbf{Ab}(O) $$ is defined by 
$$x\mapsto (x^{p}, Ab^{'}(x)) $$ 
and the map $$ \mathbf{p}:O_{\mathbf{p}}\rightarrow O_{\mathbf{p}} $$ is given by
$$(x,b)\mapsto (x^{p},p.b) $$ clearly this map factors in the category of topological spaces as 
$$ O_{\mathbf{p}}\rightarrow O^{'}\rightarrow O_{\mathbf{p}} $$
$$(x,b)\mapsto x \mapsto (x^{p},p.Ab^{'}(x)) $$
where $b=Ab^{'}(x)$.\\
Thus the abelian groups $\pi_{\ast}(O_{\mathbf{p}})$ are $p$-torsion, in particular $\tilde{\mathrm{H}}_{\ast}(\B O_{\mathbf{p}};\mathbf{Z})$ are $p$-torsion hence $\tilde{\mathrm{H}}_{\ast}(O_{\mathbf{p}};\mathbf{Z})$ are $p$-torsion ( $\B O_{\mathbf{p}}$ is simply connected, we apply the Serre spectral sequence for homology and \cite[theorem 33.68]{strom2011modern}).
At this stage we have that 
\begin{enumerate}
\item By hypothesis $R^{\delta}\rightarrow O$ is an embedding and $O$ is precofibrant, hence $R\rightarrow \Lc O$ is cofibration of topological groups, thus $\Lc O$ is Hausdorff but $\Lc O\rightarrow \Lc (O^{'})$ is a cofibration, consequently $\Lc (O^{'})$ is Hausdorff (cf. \ref{cw-complex}), in particular $O^{'}$ is Hausdorff \ref{coftgroup}. Therefore $O^{'}$ is a principal $R^{\delta}$-bundle.

\item  Since $O$ is precofibrant, $1\ast_{R^{\delta}}O$ is a cofibrant topological group, thus $\mathbf{Ab}(1\ast_{R^{\delta}}O)$ is a cofibrant topological abelian group, in particular $\mathbf{Ab}(O)/R^{\delta}$ is a Hausdorff abelian group \ref{cw-complex}. It follows that $R^{\delta}\rightarrow \mathbf{Ab}(O)$ is a closed embedding, hence $\mathbf{Ab}(O)$ is a Hausdorff principal $R^{\delta}$-bundle.

\item $O_{\mathbf{p}}$ is defined as a pullback, it is clear that $O_{\mathbf{p}}$ is a Hausdorff principal $R^{\delta}$-bundle.
\end{enumerate}
\textbf{Step 2:}\\
Applying the functor $\E := R\times_{R^{\delta}} - $ to the pullback diagram \ref{eqpullbak}, we obtain 
again a pullback diagram in the category of topological space spaces given by 
\begin{equation}\label{eqpullbak2}
\xymatrix{\E O_{\mathbf{p}}\ar[r]\ar[d]^-{inj} & \E\mathbf{Ab}(O)\ar[d]^-{\mathbf{p}}\\
\E O^{'}\ar@{->>}[r]^-{Ab^{'}} & \E\mathbf{Ab}(O)
}
\end{equation} 
where the map $\E O^{'}\rightarrow \E \mathbf{Ab}(O)$ is again a fibration. In order to demonstrate it, we proceed by step, we have proved that $\mathbf{Ab}(O)$ and $O^{'}$ are principal $R^{\delta}$-bundles. Let assume that they are trivial bundles, then up to isomorphism the pullback diagram \ref{eqpullbak} is given by 
$$\xymatrix{ L\ar[r]\ar[d] & R^{\delta}\times U\ar[d]^-{\mathbf{p}\times g}\\
R^{\delta}\times V\ar[r]^{id\times f } & R^{\delta}\times U
}$$
with $\mathbf{p}\times g: R^{\delta}\times U\rightarrow R^{\delta}\times U$ is given by $\mathbf{p}\times g(r,b)=(p.r,g(b))$, where $\mathbf{p}$ is multiplication by a prime number $p$, and 
$id \times f: R^{\delta}\times V\rightarrow R^{\delta}\times U$ is given by $id \times f (r,y)=(r,f(y)). $
Thus $$L=\{ (r,v,\frac{r}{p},u) \subset R^{\delta}\times V\times R^{\delta}\times U : f(v)=g(u) \}$$ where the action of $R^{\delta}$ on $L$ is given by $R^{\delta}\times L\rightarrow L$,$$s,(r,v,\frac{r}{p},u)\mapsto (s+r,v,\frac{s+r}{p},u).$$ It follows that 
$\E L= R\times V\times_{R\times U}R\times U=\E(R^{\delta}\times V)\times_{\E(R^{\delta}\times U)}\E(R^{\delta}\times U)$. Moreover the map $\E(R^{\delta}\times V)\rightarrow\E(R^{\delta}\times U) $ is a fibration since $R^{\delta}\times V\rightarrow R^{\delta}\times U $ is a fibration. Hence the diagram \ref{eqpullbak2} is a pullback diagram since it is a pullback diagram on local trivializations charts. Moreover the diagram \ref{eqpullbak2} is a homotopy pullback since $\E(O^{'})\rightarrow \E\mathbf{Ab}(O)$ is a fibration.

\textbf{Step 3:}\\ We observe that 
\begin{enumerate} 
\item $\E O^{'}\sim \B R^{\delta}$ (weak homotopy equivalence of spaces).
\item the action of $R^{\delta}$ on $\mathrm{H}_{\ast}(O_{\mathbf{p}};\mathbf{Z})$ is trivial (since $O_{\mathbf{p}}$ is path connected and the subgroup $R^{\delta}$ acts by multiplication on the left), moreover the action is free and proper (since $O_{\mathbf{p}} $ is a principal $R^{\delta}$-bundle) therefore by Cartan-Leray spectral sequence \cite[Theorem $8^{bis}9$]{mccleary2001user} $$E^{2}_{\ast,\ast}=\mathrm{H}_{\ast}(\B R^{\delta}; \mathrm{H}_{\ast}(O_{\mathbf{p}};\mathbf{Z}))$$ converges to 
$\mathrm{H}_{\ast}(\E O_{\mathbf{p}};\mathbf{Z})$ 
\end{enumerate}
But the homology groups $\tilde{\mathrm{H}}_{\ast}(O_{\mathbf{p}};\mathbf{Z})$ are $p$-torsion it follows by \cite[Lemma 1.1]{suslin1983thek} that $$\mathrm{H}_{\ast}(\B R^{\delta}; \mathrm{H}_{\ast}(O_{\mathbf{p}};\mathbf{Z}))=\mathrm{H}_{\ast}(\B R^{\delta};\mathbf{Z}),$$ in particular the spectral sequence collapses at the second page, hence 
$$\mathrm{H}_{\ast}(\E O_{\mathbf{p}};\mathbf{Z})=\mathrm{H}_{\ast}(\B R^{\delta};\mathbf{Z}).$$
The homology groups $\tilde{\mathrm{H}}_{\ast}(\B R^{\delta};\mathbf{Z})$ are rational vector spaces, thus by universal coefficients theorem $\tilde{\mathrm{H}}_{\ast}(\E O_{\mathbf{p}};k)=0$ for any finite field $k$. Since $\E \mathbf{Ab}(O)$ is a topological group (because $ \mathbf{Ab}(O)$ is abelian \ref{cartesian}), we have the following fiber sequence

$$\xymatrix{
\E O_{\mathbf{p}}\ar[r] & \E O^{'} \times \E \mathbf{Ab}(O)\ar[r]^-{f} & \E \mathbf{Ab}(O)
}$$
where $f(a,b)= Ab(a)[\mathbf{p}(b)]^{-1}.$
But $\tilde{\mathrm{H}}_{\ast}(\E O_{\mathbf{p}};k)=0$ for any finite field $k$, applying the Serre spectral sequence to the homotopy fiber we obtain that
$$\mathbf{p}_{\ast}: \mathrm{H}_{\ast}(\E \mathbf{Ab}(O);k)\rightarrow \mathrm{H}_{\ast}(\E \mathbf{Ab}(O);k)$$ is an isomorphism for any finite field $k$. But the $p$-power map  $$\mathbf{p}_{\ast}: \mathrm{H}_{\ast}(\E \mathbf{Ab}(O);\mathbf{Q})\rightarrow \mathrm{H}_{\ast}(\E \mathbf{Ab}(O);\mathbf{Q})$$ induces an isomorphism hence 
$$\mathbf{p}_{\ast}: \mathrm{H}_{\ast}(\E \mathbf{Ab}(O);\mathbf{Z})\rightarrow \mathrm{H}_{\ast}(\E \mathbf{Ab}(O);\mathbf{Z})$$  is an isomorphism and therefore $\mathbf{p}:\E \mathbf{Ab}(O)\rightarrow \E \mathbf{Ab}(O)$ is a weak homotopy equivalence in particular $\E \mathbf{Ab}(O)$ is a rational topological group in the sense of \ref{rationality0}. We have noticed that $R^{\delta}\rightarrow \mathbf{Ab}(O)$ is a closed embedding of topological groups, therefore we have a fiber sequence $R^{\delta}\rightarrow R\times  \mathbf{Ab}(O)\rightarrow \E \mathbf{Ab}(O)$. We conclude that $\mathbf{Ab}(O)$ is a rational topological group, finally by lemma \ref{rationality} we conclude that $\Lc(O)$ is rational (all homotopy groups are rational vector spaces).  
\end{proof}
\begin{theorem}\label{crucial}
Suppose that $O$ is an object of $R^{\delta}/\Gr$ such that 
\begin{enumerate}
\item $O$ is contractible as a space. 
\item  $R^{\delta}\rightarrow O$, $R^{\delta}\rightarrow \mathbf{Ab}(O)$ are embeddings of topological groups.
\item  $\Lc(O)$ is a retract of a CW-complex.
\end{enumerate}
Then the natural map $$\mathrm{H}_{\ast}(O;k)\rightarrow \mathrm{H}_{\ast}(\Lc O;k)$$ is an isomorphism for any finite field $k$. 
\end{theorem}
\begin{proof}
The map $R\rightarrow \Lc(O)$ can be factored as a cofibration followed by a trivial fibration
$$\xymatrix{R\ar@{>->}[r] & D \ar@{->>}[r]^-{\sim} & \Lc(O) } $$
in the category $R/\Gr$. In particular $R\rightarrow D$ is a cofibrant object in $R/\Gr$ and $\Lc D=D$ (here we see the composition $R^\delta\rightarrow R\rightarrow D$ as an object of $R^{\delta}/\Gr$). By \ref{coftgroup} $D$ is Hausdorff. Consider the pullback diagram in the category $R^{\delta}/\Gr$ given by 

$$ \xymatrix{D^{'}\ar@{->>}[r]^-{\sim} \ar[d] & O\ar[d]\\
D\ar@{->>}[r]^-{\sim} & \Lc O}$$
Clearly $D^{'}$ is a weakly contractible Hausdorff topological group and $R^{\delta}\rightarrow D^{'}$, $R^{\delta}\rightarrow \mathbf{Ab}(D^{'})$ are embeddings. Since $\Lc (O)$ is a retract of CW-complexes it follows that 

$$ A\rightarrow D\rightarrow \Lc (O) $$
(where $A=\mathbf{ker}(D\rightarrow \Lc O)$) is a short exact sequence of topological groups and a trivial $A$-principal bundle (since $\Lc O$ is a retract of a CW-complex, we have a section $\Lc (O)\rightarrow D$). Notice that as a space $A$ is contractible. By construction $A\rightarrow D^{'}\rightarrow O$ is also a short exact sequence of topological groups and in the same time a trivial $A$-principal bundle. Thus we obtain a commutative diagram of (topologically split) short exact sequences of topological groups (where \textbf{all vertical maps are homomorphisms of topological groups and identities of underlying groups})
$$ \xymatrix{A\ar[r]\ar@{=}[d]& D^{'}\ar@{->>}[r]^-{\sim}\ar[d] & O\ar[d]\\
A\ar@{=}[d]\ar[r]& \Lc D^{'}\ar[r]\ar[d] & \Lc O\ar@{=}[d]\\
A\ar[r]& D\ar@{->>}[r]^-{\sim} & \Lc O}$$
\textbf{Here is a precise explanation of why these are topologically split short exact sequences of topological groups:}\textit{ 
since $\Lc O$ is a retract of a CW-complex and $q:D\rightarrow \Lc O$ is a trivial fibration, there exists a continuous section $s: \Lc O\rightarrow D$ such that $q\circ s=id$. It follows that $A\rightarrow D\rightarrow \Lc O$ is a short exact sequence of topological groups. By defintion $D^{'}= D\times_{\Lc O} O$, therefore $s: O\rightarrow D^{'}$ is a section to the map $q: D^{'}\rightarrow O$, hence  $A\rightarrow D^{'}\rightarrow O$ is a short exact sequence of topological groups, in particular $q:D^{'}\rightarrow  O$ is a quotient map. The map $t: D\rightarrow A\times \Lc O$ given by $d\mapsto (d.[s(q(d))]^{-1},q(d))$ is a homeomorphism, and similarly the map $t: D^{'}\rightarrow A\times O$ given by $d\mapsto (d.[s(q(d))]^{-1},q(d))$ is a homeomorphism.  Therefore the map $r: \Lc D^{'}\rightarrow A$ given by  $d\mapsto d.[s(q(d))]^{-1}$ is a continuous map. It follows that the map $\Lc D^{'}\rightarrow \Lc D^{'}$  given by $d\mapsto s(q(d))$ is continuous. Since $\Lc$ is a functor, $A\rightarrow \Lc D^{'}\rightarrow \Lc O$ is a short exact sequence of topological groups. Thus $q: \Lc D^{'}\rightarrow \Lc O$ is a quotient map, therefore $s:\Lc O\rightarrow \Lc D^{'}$ is a continuous map, in particular  the map $t: \Lc D^{'}\rightarrow A\times \Lc O$ given by $d\mapsto (d.[s(q(d))]^{-1},q(d))$ is a homeomorphism. } \\\\
It implies that the homomorphism of groups $\Lc D^{'}\rightarrow D$ is a homeomorphsism hence $\Lc D^{'}=D$, thus $D^{'}$ is precofibrant in the sense of \ref{qcofibrant} since $D$ is cofibrant in $R/\Gr$ by construction. By hypothesis $D^{'}$ is contractible, we conclude, by \ref{contractible}, that $D$ is rational hence $\Lc O$ is rational. In particular the natural map $$\mathrm{H}_{\ast}(O;k)\rightarrow \mathrm{H}_{\ast}(\Lc O;k)$$ is an isomorphism for any finite field $k$. 
\end{proof}

\begin{theorem}\label{crucial+}
Let $i:R^{\delta}\rightarrow G$ be an object in $R^{\delta}/\Gr$  such that $\Lc(G)$ is a retract of a CW-complex as a space, then the natural map $$\mathrm{H}_{\ast}(\B G;k)\rightarrow \mathrm{H}_{\ast}(\B \Lc G;k)$$ is an isomorphism for any finite field $k$. 
\end{theorem}
\begin{proof}
We consider the topological group $H=G\times R^{\delta}$ and 
$i:R^{\delta}\rightarrow H$ is given by the composition $(i\times id)\circ diag$

$$\xymatrix{R^{\delta}\ar[r]^-{diag} & R^{\delta}\times R^{\delta} \ar[r]^-{i\times id} & G \times R^{\delta} }$$ 

The map $H\rightarrow 1$ can be factored in  $R^{\delta}/\Gr$ as $H\rightarrow O\rightarrow 1$ where the first map is a cofibration and the second map is a trivial fibration. We have the following  pushout diagram in the category of topological groups (in fact a homotopy pushout diagram since the category of topological groups is left proper) given by 

$$\xymatrix{H\ar[r]\ar@{>->}[d] & \Lc H\ar[d]\\
O\ar[r] & \Lc O \cocartesien}$$
The upper horizontal map is the identity map in the category of groups, therefore 
\begin{equation} \label{tot}
 \xymatrix{\B H\ar[r]\ar[d] & \B\Lc H\ar[d]\\
\B O\ar[r] & \B\Lc O}
\end{equation}
is a homotopy pushout in the category of topological spaces \ref{B-pushout}. Clearly $O$ verifies the following properties:
\begin{enumerate}
\item $O$ is contractible.
\item $R^{\delta}\rightarrow O$, $R^{\delta}\rightarrow \mathbf{Ab}(O)$ are embeddings of topological groups.
\item $H\rightarrow O$ is a cofibration of topological groups hence $\Lc (H)\rightarrow \Lc (O)$ is cofibration of topological groups and since $\Lc(H)= \Lc(G)\times R$ is a retract of a CW-complex, then by corollary \ref{coftgroup}, $\Lc (O)$ is a retract of a $CW$-complex.
\end{enumerate}
Applying theorem \ref{crucial}, we have that 
$$\mathrm{H}_{\ast}(O;k)\rightarrow \mathrm{H}_{\ast}(\Lc O;k)$$ is an isomorphism for any finite field $k$, it follows that 
$$\mathrm{H}_{\ast}(\B O;k)\rightarrow \mathrm{H}_{\ast}(\B \Lc O;k)$$ is an isomorphism for any finite field $k$, in particular, by Meyer-Vietoris long exact sequence applied to the homotopy pushout diagram \ref{tot}, the natural map $$\mathrm{H}_{\ast}(\B H;k)\rightarrow \mathrm{H}_{\ast}(\B \Lc H;k)$$  is an isomorphism for any finite field $k$. Finally the result follows since $\Lc H= \Lc G\times R$, in particular  the natural map $$\mathrm{H}_{\ast}(\B G;k)\rightarrow \mathrm{H}_{\ast}(\B \Lc G;k)$$  is an isomorphism for any finite field $k$.
\end{proof}

\begin{remark}
The third hypothesis  of theorem \ref{crucial} seems to be restrictive, but it is sufficient for proving the Friedlander-Milnor conjectures. Meanwhile one can hope for a more general result without assuming the third hypothesis (CW-complex). At this stage, we are unable to provide an evidence for more global result except the results of section \ref{section4}.   
\end{remark}
\section{Lie Groups}\label{section6}
 \begin{warning}
Given a topological group or a topological monoid $G$ with two different topologies $\tau_{2}<\tau_{1}$ where $\tau_{1}$ is finer then $\tau_{2}$
we will abusively denote the evident continuous map $G^{\tau_{1}}\rightarrow  G^{\tau_{2}}$ (identity on the underlying algebraic structures) by $id$. The reader should be aware that $id$ is bijective continuous map but not an homemorphism of underlying topological spaces.  
\end{warning}
\begin{proposition}\label{diffeo}
Let $G$ be a connected real simple Lie group of dimension $n$ (\cite{tits1967tabellen}) and let $SO(2)=T\rightarrow G$ be an injective homomorphism of topological groups (which has to exist), then 
there exists $n$ elements $g_{1},\dots, g_{n}\in G$ such that the smooth map 
$$\mu:T_{1}\times\dots \times T_{n}\rightarrow G $$
$$ (x_{1},\dots, x_{n})\mapsto x_{1}x_{2}\dots x_{n} $$
is a local diffeomorphism in a neighbourhood of $(e_{1},\dots e_{n})$, where  $T_{i}=g_{i}Tg_{i}^{-1}$ and $e_{i}$ is the neutral (unit) element of $T_{i}$. 
\end{proposition}
 \begin{proof}
We start with some standard notation, the (real) Lie algebra associated to the real Lie group $G$ is denoted by $\mathfrak{g}$, recall that as a vector space the Lie algebra $\mathfrak{g}$ is identified to the tangent space $\mathbf{T}_{e}G$. Since $G$ is a connected simple Lie group, it follows that the center of $G$ is trivial. Hence, the adjoint representation $\rho: G\rightarrow \mathrm{GL}(\mathfrak{g})$ is continuous, injective and closed homomorphism of Lie groups. Therefore without loosing generality, we can suppose that $G$ is a closed subgroup of a general linear group $\mathrm{GL}_{n}(\mathbf{R})$. In particular, elements of $G$ are matrices and the multiplication is multiplication of matrices. 

 For any element $g\in G$ we define the inner homomorphism of groups $c_{g}:G\rightarrow G$ by $x\mapsto gxg^{-1}$. The derivative of the morphism $c_{g}$ at the neutral element $e$ is denoted by $Ad_{g}:= \mathbf{d}c_{g}(e)$. Let $v\in \mathbf{T}_{e}G$, where $v$ is a non trivial tangent vector to $T\subset G$ at $e$. We consider the set 
 $$S =\{\mathbf{d}c_{g}(e)(v)=gvg^{-1}| g\in G\}.$$
Let $ V= Span\{S\}$ (i.e. the set $S$ generates $V$ as a vector space), clearly $G$ acts by conjugation on $V$ (which is not trivial since $v\neq 0$ ), hence we have a sub-representation $G\rightarrow \mathrm{GL}(V)$ of the adjoint representation $ad: G\rightarrow\mathrm{GL}(\mathfrak{g})$, but $G$ is simple therefore the adjoint representation of $G$ is irreducible. It implies that $V=\mathfrak{g}$. Since $\mathfrak{g}$ is a vector space of dimension $n$, there exists a finite family $\{g_{1},\dots , g_{n}\}$ such that $\{Ad_{g_{1}}(v),\dots , Ad_{g_{n}}(v)\}$ form a basis for the vector space $\mathfrak{g}$. In particular 
 $$\mathfrak{g}=  V_{1}\oplus V_{2}\oplus \dots \oplus V_{n}$$
 where $V_{i}$ is the one dimensional vector space generated by $Ad_{g_{i}}(v)$.
We define $T_{i}$ by $g_{i}Tg_{i}^{-1}$ and consider the smooth morphism of smooth manifolds:
$$\mu:T_{1}\times\dots \times T_{n}\rightarrow G $$ 
where each $T_{i}$ is embedded naturally in $G$ and each $e_{i}$ will denote the neutral element of $T_{i}$, then by \cite{lee2009manifolds} the linear morphism (which we describe up to linear isomorphism)$$ \mathbf{d}\mu(e_{1},\dots ,e_{n}):\mathbf{T}_{(e_{1},\dots ,e_{n})}T_{1}\times \dots \times T_{n} \cong V_{1}\oplus V_{2}\oplus \dots \oplus V_{n}\longrightarrow \mathfrak{g}=\mathbf{T}_{e}G$$
 sends any $(v_{1},\dots  ,v_{n})\in \mathbf{T}_{(e_{1},\dots ,e_{n})}T_{1}\times \dots \times T_{n}$ to
 $$\mathbf{d}\mu_{(e_{1},\dots e_{n})}(v_{1},\dots, v_{n})= v_{1}+ \dots + v_{n}. $$
 In particular $\mathbf{d}\mu_{(e_{1},\dots e_{n})}$ has rank equal to $n$, therefore by inverse function theorem there exists an open neighbourhood $U$ of $(e_{1},\dots e_{n})$ in $T_{1}\times \dots \times T_{n}$ (equipped with the product topology) such that $\mu:U\rightarrow \mu(U)$ is a diffeomorphism \cite[Theorem C.1]{lee2009manifolds}. 
 
 \end{proof}

\begin{theorem}\label{pushout}
Let $G$ be a connected simple Lie group of dimension $n$ (hence the maximal torus is non trivial), then the commutative square 
$$\xymatrix{ R^{\delta}\ar[r]\ar[d]_{id} & G^{\delta}\ar[d]^{id}\\
R\ar[r] & G
}$$
is a pushout in the category $\Gr$ i.e. $\Lc (G^{\delta})= G$.
\end{theorem}
\begin{proof}
First of all, the precedent commutative square is a pushout in the category $\Grp$ ones we forget about topology, it remains to prove that $G$ (with multiplication $\mu$) with its standard topology of Lie group is the pushout of the diagram in  $\Gr$. Since $G$ is simple there is an embedding of topological groups $T\rightarrow G $. Let suppose that the dimension of $G$ as a Lie group is $n$ and lets suppose that the pushout $G^{\tau}$, 
$$\xymatrix{ R^{\delta}\ar[r]\ar[d]_{id} & T^{\delta} \ar[r]\ar[d]_{id} & G^{\delta}\ar[d] \ar[rdd]&\\
R\ar[r]\ar[drrr] &  T\cocartesien \ar[r]\ar[drr] & G^{\tau}\cocartesien \ar@{.>}[dr]^-{\exists !}_-{u} & \\
& & & G
}$$ 
is endowed with Hausdorff topology $\tau$ (since $G^{\tau}\rightarrow G$ is continuous and bijective and $G$ is Hausdorff) different from the one making it as Lie group. By hypothesis we have a continuous map $T\rightarrow G^{\tau}$, which allows us to define a map $\mu^{'} : T_{1}\times \dots \times T_{n}\rightarrow G^{\tau}$ by $\mu^{'}(x_{1},\dots, x_{n})=x_{1}x_{2}\dots x_{n}\in G^{\tau}$, where the product $T_{1}\times \dots \times T_{n}$ is endowed with the product topology ($n$-dimensional torus) and $T_{i}=g_{i}Tg_{i}^{-1}$ defined as in \ref{diffeo}. The map $\mu^{'}$ is obviously continuous. On another hand, the composition $\mu=u\circ\mu^{'} $ given by
$$\xymatrix{T_{1}\times \dots \times T_{n}\ar[r]^-{\mu^{'}} & G^{\tau}\ar[r]^-{u} & G}$$ 
is continuous. It follows that there exists a small neighbourhood $U$ of $(e_{1},e_{2},\dots, e_{n})\in U\subset  T_{1}\times\dots \times T_{n}$ such that $\mu_{|U} :U\rightarrow \mu(U)$ is a diffeomorphism by \ref{diffeo}. Therefore $\mu^{'}_{|U}:U\rightarrow \mu^{'}(U)\subset G^{\tau} $ is a continuous bijection. 
But $T_{1}\times \dots \times T_{n}$ is compact and $G^{\tau}$ is Hausdorff, it implies that $\mu^{'}_{|U}$ takes a closed subset to a closed subset, which means that  $\mu^{'}_{|U}:U\rightarrow \mu^{'}(U)\subset G^{\tau} $ is a homeomorphism. We conclude that $u:G^{\tau}\rightarrow G$ is a bijective continuous map and a local homeomorphism between topological groups and hence a homeomorphism i.e., an isomorphism in the category $\Gr$.

\end{proof}
\section{Proof of the Friedlander-Milnor conjectures (F-M)}\label{section7}

\begin{theorem}\label{main}
The Milnor conjecture is true, consequently the Friedlander conjecture is true for complex algebraic Lie groups.
\end{theorem}
\begin{proof}
Let $G$ be a Lie group and $G^{\delta}$ be the same group with discrete topology, in order to prove the Milnor conjecture it is sufficient to prove it for simply connected simple Lie groups \cite[Lemma 4]{milnor1983homology}. Let $G$ be a connected simple Lie group then according to Theorem \ref{pushout} it follows that $\Lc(G^{\delta})=G$ in $\Gr$ ($G$ is a CW-complex). We consider the diagonal embedding $R^{\delta}\rightarrow G^{\delta}\times R^{\delta}$, therefore $\Lc(G^{\delta}\times R^{\delta})=\Lc(G^{\delta})\times R$. Hence by applying  Theorem \ref{crucial+}, we conclude that
$$\mathrm{H}_{\ast}(\B G^{\delta}; \mathbf{F}_{p})=\mathrm{H}_{\ast}(\B G^{\delta}\times \B R^{\delta}; \mathbf{F}_{p})\rightarrow  \mathrm{H}_{\ast}(\B \Lc(G^{\delta})\times \B R; \mathbf{F}_{p})= \mathrm{H}_{\ast}(\B G; \mathbf{F}_{p})$$
is an isomorphism in homology with coefficients in $\mathbf{F}_{p}$ for any prime number $p$. The same result hold for Mod-$p$ cohomology and hence the Friedlander conjecture over $\mathbf{C}$ follows in the case of complex algebraic Lie groups by remark \ref{remark0}. 
\end{proof}
\begin{remark}
Recall that the Milnor conjecture is false with rational coefficients \cite[Lemma 8]{milnor1983homology} i.e., the induced map $\mathrm{H}_{\ast}(\mathbf{B}G^{\delta};\mathbf{Q})\rightarrow \mathrm{H}_{\ast}(\mathbf{B}G;\mathbf{Q})$ is not an isomorphism in general. 
\end{remark}
\begin{theorem}(Stable Milnor conjecture)\label{stable}
 Let $p$ be any prime number and $\mathbf{k}$ the filed of real numbers $\mathbf{R}$ or the field of complex numbers $\mathbf{C}$. We define $\B GL_{\infty}(\mathbf{k})$ as $ colim_{n}\B GL_{n}(\mathbf{k})$, then the evident map $$\Homology_{\ast}(\B GL_{\infty}(\mathbf{k})^{\delta};\mathbf{F}_{p})\rightarrow \Homology_{\ast}(\B GL_{\infty}(\mathbf{k});\mathbf{F}_{p})$$ is an isomorphism. The same result holds if we replace $GL_{\infty}$ by the special linear group $SL_{\infty}$, by the orthogonal group $O_{\infty}$, by the special orthogonal group $SO_{\infty}$ or by the unitary group $U_{\infty}$, by the special unitary group $SU_{\infty}$ or by the symplectic group $Sp_{\infty}$.
\end{theorem}
\begin{proof}
The proof is the same for all cited groups. We will do it for $GL_{\infty}$. We consider the following commutative diagram 
$$\xymatrix{ \B GL_{1}(\mathbf{k})^{\delta}\ar[r] \ar[d] & \B GL_{2}(\mathbf{k})^{\delta}\ar[r] \ar[d] & \dots \ar[r] \ar[d] & \B GL_{n}(\mathbf{k})^{\delta}\ar[r] \ar[d] &\dots\\
\B GL_{1}(\mathbf{k})\ar[r] & \B GL_{2}(\mathbf{k})\ar[r] & \dots \ar[r] & \B GL_{n}(\mathbf{k})\ar[r]  & \dots
} $$
The colimit of the upper diagram is by definition $\B GL_{\infty}(\mathbf{k})^{\delta}$ which is also a homotopy colimit in the category $\Top$. In the same way the colimit of the lower diagram is by definition $\B GL_{\infty}(\mathbf{k})$ which is also a homotopy colimit in the category $\Top$. By theorem \ref{main} the induced map $$C_{\ast}(\B GL_{n}(\mathbf{k})^{\delta};\mathbf{F}_{p})\rightarrow C_{\ast}(\B GL_{n}(\mathbf{k});\mathbf{F}_{p})$$ is a quasi isomorphism of chain complexes for any $n\in \mathbf{N}$. Therefore taking the homotopy colimit in the model category of chain complexes, the natural map 
$$\hocolim_{n} C_{\ast}(\B GL_{n}(\mathbf{k})^{\delta};\mathbf{F}_{p})\rightarrow \hocolim_{n} C_{\ast}(\B GL_{n}(\mathbf{k});\mathbf{F}_{p}) $$ 
is a quasi isomorphism ($\sim$) of chain complexes. But $\hocolim_{n} C_{\ast}(\B GL_{n}(\mathbf{k})^{\delta};\mathbf{F}_{p})\sim C_{\ast}(\B GL_{\infty}(\mathbf{k})^{\delta};\mathbf{F}_{p})$ and $\hocolim_{n} C_{\ast}(\B GL_{n}(\mathbf{k});\mathbf{F}_{p})\sim C_{\ast}(\B GL_{\infty}(\mathbf{k});\mathbf{F}_{p})$. In particular the natural map of chain complexes $$ C_{\ast}(\B GL_{\infty}(\mathbf{k})^{\delta};\mathbf{F}_{p})\rightarrow C_{\ast}(\B GL_{\infty}(\mathbf{k});\mathbf{F}_{p}) $$ is a quasi isomorphism. We conclude that the induced map
$$ \Homology_{\ast}(\B GL_{\infty}(\mathbf{k})^{\delta};\mathbf{F}_{p})\rightarrow \Homology_{\ast}(\B GL_{\infty}(\mathbf{k});\mathbf{F}_{p}) $$
is an isomorphism and therefore $\Homology^{\ast}(\B GL_{\infty}(\mathbf{k});\mathbf{F}_{p})\rightarrow \Homology^{\ast}(\B GL_{\infty}(\mathbf{k})^{\delta};\mathbf{F}_{p}) $ is an isomorphism.

\end{proof}

Some applications of Theorem \ref{stable} overlap with some results about the computational aspects of the algebraic $K$-theory of local fields and algebraically closed fields of characteristic 0 described in \cite{suslin1984k} and \cite{suslin1983thek}.


\bibliographystyle{plain}
\bibliography{spb8}
\end{document}